\documentclass[12pt]{article}
\usepackage{a4wide}
\usepackage{amssymb}

\usepackage{latexsym}
\usepackage{graphics}
\usepackage{amsfonts}
\usepackage{color}

\newcommand{\ben}{\begin{enumerate}}
\newcommand{\een}{\end{enumerate}}
\newcommand{\ble}{\begin{lem}}
\newcommand{\ele}{\end{lem}}
\newcommand{\bth}{\begin{thm}}
\renewcommand{\eth}{\end{thm}}
\newcommand{\bpr}{\begin{prop}}
\newcommand{\epr}{\end{prop}}
\newcommand{\bco}{\begin{cor}}
\newcommand{\eco}{\end{cor}}
\newcommand{\bcon}{\begin{conj}}
\newcommand{\econ}{\end{conj}}

\newcommand{\bde}{\begin{defn}}
\newcommand{\ede}{\end{defn}}

\newcommand{\bex}{\begin{exa}}
\newcommand{\eex}{\end{exa}}
\newcommand{\barr}{\begin{array}}
\newcommand{\earr}{\end{array}}
\newcommand{\btab}{\begin{tabular}}
\newcommand{\etab}{\end{tabular}}
\newcommand{\beq}{\begin{equation}}
\newcommand{\eeq}{\end{equation}}
\newcommand{\bea}{\begin{eqnarray*}}
\newcommand{\eea}{\end{eqnarray*}}
\newcommand{\bce}{\begin{center}}
\newcommand{\ece}{\end{center}}
\newcommand{\bpi}{\begin{picture}}
\newcommand{\epi}{\end{picture}}
\newcommand{\bfi}{\begin{figure} \begin{center}}
\newcommand{\efi}{\end{center} \end{figure}}

\newcommand{\bsl}{\begin{slide}{}}
\newcommand{\esl}{\end{slide}}

\newcommand{\pf}{\noindent{\bf Proof}\hspace{7pt}}
\newcommand{\qed}{\rule{1ex}{1ex}}

\newcommand{\hso}[1]{\hspace{-1pt}}

\newcommand{\into}{\hookrightarrow}

\def\<{\langle}
\def\>{\rangle}

\newcommand{\cA}{{\cal A}}

\newcommand{\cC}{{\cal C}}

\newcommand{\Aut}{\mathop{\rm Aut}\nolimits}

\newcommand{\ad}{\mathop{\rm ad}\nolimits}

\renewcommand{\bar}{\overline}

\newcommand{\id}{\mathop{\rm id}\nolimits}
\newcommand{\im}{\mathop{\rm im}\nolimits}

\newcommand{\orG}{\overrightarrow{\Gamma}}

\newcommand{\orE}{\overrightarrow{E}}

\newcommand{\dfn}{\em}

\newcommand{\after}{\mathbin{ \circ }}

\renewcommand{\hat}{\widehat}

\newcommand{\mn}{\ \medskip \newline }
\newtheorem{thm}{Theorem}[section]
\newtheorem{prop}[thm]{Proposition}
\newtheorem{cor}[thm]{Corollary}
\newtheorem{lem}[thm]{Lemma}
\newtheorem{conj}[thm]{Conjecture}
\newtheorem{exa}[thm]{Example}

\newtheorem{defn}[thm]{Definition}

\newtheorem{mainthm}{Theorem}

\renewcommand{\qed}{\hfill $\square$}

\newcommand{\Note}{\ \newline\refstepcounter{thm}\noindent{\bf Note \thethm}\quad}

\newcounter{romanlistctr}
{\end{list}}%

 \makeatletter
 \@addtoreset{equation}{section}
 \makeatother

 \makeatletter
 \def\section{\@startsection {section}{1}{\z@}{-1.5ex plus -.5ex
 minus -.2ex}{1ex plus .2ex}{\large\bf}}

 \makeatletter
 \def\subsection{\@startsection {subsection}{1}{\z@}{-1.5ex plus -.5ex
 minus -.2ex}{1ex plus .2ex}{\bf}}

\def\<{\langle}
\def\>{\rangle}

\title{Bass-Serre theory and counting rank two amalgams.}
\author{Rieuwert J. Blok $^1$ and Corneliu Hoffman$^2$\footnote{This paper was written during a Research in Pairs visit to the Mathematisches Forschungsinstitut Oberwolfach. We would like to thank the institute for their hospitality and generous support.}\\
$^1$Department of Mathematics and Statistics
\\ Bowling Green State
University \\ Bowling Green, OH 43403-1874\\[5pt]
$^2$ School of Mathematics \\
University of Birmingham \\
Edgbaston, B15 2TT, United Kingdom \\
{\tt blokr@member.ams.org}\\[5pt]
}

\date{Submitted: October 29, 2009\\[1in]
    \begin{flushleft}
    \end{flushleft}
       }
 \makeindex

\begin{document}
\pagestyle{empty}
 \maketitle
\section{Introduction}
An amalgam of groups can be viewed as a Sudoku game inside a group. You are given a set of subgroups and their intersections and you need to decide what the largest group containing such a structure can be. This approach is very useful for example in the classification of finite simple groups. More precisely,  induction and local analysis provides a set of subgroups of the minimal counterexample and then amalgam type results such as the Curtis-Tits and Phan theorems show that the group is known after all. 

Recently the Curtis-Tits-Phan theory has seen a lot of activity. Similar results have been obtained for other finite groups of Lie type and more generally some   subgroups of  Kac-Moody groups (see \cite{BenGraHofShp03,Gra} for details). All these results are similar in nature and prove that a certain group is the universal completion of a rank 2 amalgam of its subgroups. 

This leaves open the question of whether just the structure of the subgroups involved determines the group. Most approaches to this problem use induction together with a lemma by Goldschmidt (see Corollary~\ref{lem:Goldschmidt}) that describes the isomorphism classes of amalgams of two groups in terms of double coset enumeration. Our results will generalise this result.

 In a recent work \cite{BloHof09a} we used Bass-Serre theory of graphs of groups to classify all possible amalgams of Curtis-Tits shape with a given diagram. This note describes the  method for general rank two amalgams. In Section~\ref{sec:Amalgams and graphs of groups} we introduce the main notions, prove Theorem~\ref{thm:pointings and amalgams}  and obtain Goldschmidt's lemma as a particular case. Moreover, we specialise to the case of a triangle to get a very concrete application to amalgams coming from rank three geometries. In Section~\ref{sec:fundamental group} we specialise to what we shall call rigid amalgams and discuss fundamental groups and their significance for classifying such amalgams.

\section{Amalgams and graphs of groups}\label{sec:Amalgams and graphs of groups}

In this paper we fix a graph $\Gamma=(I,E)$ with vertex set $I$ and edge set $E$, where multiple edges and loops are permitted.

\bde\label{dfn:amalgam}
An {\em amalgam} over the  graph $\Gamma=(I,E)$ is a collection $\cA=\{G_i,G_e,\varphi_{i,e}\mid i\in I,i\in e\in E\}$, where $G_i$ and $G_e$ are groups and, for each pair $(i,e)$ such that the vertex $i$ lies on the edge $e$, we have a monomorphism $\varphi_{i,e}\colon G_i\into G_e$, called an {\em inclusion map}. 
A {\em completion} of $\cA$ is a group $G$ together with a collection  $\phi=\{\phi_e,\phi_i\mid i\in I, e\in E\}$ of homomorphisms $\phi_e\colon G_e\to G$, $\phi_i\colon  G_i\to G$ such that for any $i\in e$ we have $\phi_e\after\varphi_{i,e}=\phi_i$.
For simplicity we shall write $ \bar G_{i,e}=\varphi_{e}(G_i)\le G_{e}$.
The amalgam $\cA$ is {\em non-collapsing} if it has a non-trivial completion.
A completion $(\hat{G},\hat{\phi})$ is called {\em universal} if for any completion $(G,\phi)$ there is a  (necessarily unique) surjective group homomorphism $\pi\colon \hat{G}\to G$ such that $\phi=\pi\after\hat{\phi}$. 
\ede

\mn
Following \cite{Bas93,Se1980} we define a directed graph $\orG=(I,\orE)$ where for each edge $e\in E$ we introduce directed edges $e$ and $\bar e$ in $\orE$.  Moreover we denote by $d_0(e)$ the starting node of the oriented edge $e$ and by $d_1(e)$ the end vertex of $e$. Thus $d_0(e)=d_1(\bar e)$.

From now on, we shall index every amalgam $\cA=\{G_i,G_e,\varphi_{i,e}\mid i\in I, e\in E,\mbox{ where }i\in e\}$ by vertices $i\in I$ and {\em directed} edges $e\in \orE$ and set $G_e=G_{\bar{e}}$ and $\varphi_e=\varphi_{d_0(e),e}$. Then, $\varphi_{\bar{e}}=\varphi_{d_1(e),e}$.

\bde\label{dfn:CT hom}
A {\dfn homomorphism} between the amalgams 
$\cA^1(\orG)=\{G_i^1,G_e^1,\varphi_{e}^1\mid i\in I, e\in \orE\}$ and 
$\cA^2(\orG)=\{G_i^2,G_e^2,\varphi_{e}^2\mid i\in I, e\in \orE\}$ is a
map $\phi=\{ \phi_i, \phi_{e}\mid i\in I^1, e \in E^1\}$ where $\phi_i\colon G_i^1\to  G_{i}^2$ and $\phi_{e}  \colon G_{e}^1\to G_{e}^2$ are group homomorphisms  such that $$\phi_{e}\after \varphi_{e}^1 =\varphi_{e}^2\after  \phi_{d_0(e)}.$$
We call $\phi$ an {\dfn isomorphism} of amalgams if $\phi_i$ and $\phi_{e}$ are bijective for all $i,j\in I$, and $\phi^{-1}$ is a homomorphism of amalgams. 
 \ede
 
\bde\label{dfn:type of amalgam}
Consider an amalgam $\cA_0=\{G_i,G_{e},\psi_{e}\mid i\in I, e\in \orE\}$ over $\orG=(I,\orE)$. An {\em amalgam of type $\cA_{0}$} is an amalgam over $\orG$ with the same groups $G_i, G_{e}$ and where, for every  $e\in \orE$, the inclusion map $\varphi_{e}$ is possibly different from $\psi_{e}$ except that it has the same image $\bar G_{{i,e}}$.
\ede

\mn
The aim of this note is to describe the isomorphism classes of all amalgams that have the same type as $\cA_{0}$. For the rest of the paper, $\cA_0$ will be a fixed amalgam over $\orG=(I,\orE)$.
We shall therefore simply write $\cA_0=\{G_i,G_e,\psi_{e}\}$ tacitly understanding that $G_i$ is taken over all $i\in I$, that $G_e$, is taken over all $e\in E$, and $\psi_{e}$ is taken over all edges $e\in \orE$ and adopt a similar shorthand for amalgams $\cA$ of type $\cA_0$.

\bde  \label{def:graph of groups} A {\em graph of groups} is a pair $(\cC, \orG)$ where $\orG$ is a graph as above and $\cC$ associates to each $i\in I$ a group $A_i$ and to each directed edge $e\in \orE$ a group $A_e=A_{\bar e}$.
Moreover, for each vertex $i$ with $d_0(e)=i$ we have a homomorphism $\alpha_{e}\colon A_{e}\to A_i$.

Since we shall work with a fixed graph $\orG=(I,\orE)$, if we want to specify the groups of $(\cC,\orG)$ we shall write  
 $\cC=\{A_i,A_e,\alpha_e\}$ rather than $\cC=\{A_i,A_e,\alpha_e\mid i\in I, e\in \orE\}$.
 \ede

\mn
Note that Definition~\ref{def:graph of groups} is a generalisation of the concept from \cite{Bas93} as we do not require the $\alpha_e$'s to be monomorphisms.

\bde\label{dfn:inner morphism}
Given graphs of groups $(\cC^{(k)},\orG)$ for $k=1,2$, an {\em inner morphism} is a collection $\phi=\{\phi_i,\phi_{e}\mid i\in I, e\in \orE\}$ of group homomorphisms $\phi_i\colon A_i^{(1)}\to A^{(2)}_{i}$
  and  
$\phi_{e}\colon A_{e}^{(1)}\to A^{(2)}_{e}$
  so that for each $e\in \orE$ with $d_0(e)=i$ there exists an element
   $\delta_{e}\in A_{i}^{(2)}$ so that
   $$\phi_{i}\after \alpha_{e}=\ad(\delta_{e}^{-1})\after\alpha_{e}\after \phi_{e}.$$
Here $\ad(x)(y)=x^{-1}yx$.
\ede

\bde\label{dfn:a pointing of a graph of groups}
Let $(\cC_0,\orG)$ be a graph of groups, a {\em pointing} is a 
pair $((\cC,\orG),\delta)$, where $\delta=\{\delta_{e}\mid e\in \orE\}$ is a collection  of elements $\delta_{e}\in A_{d_0(e)}$ and  $(\cC,\orG)$ is a graph of groups obtained from $(\cC_0,\orG)$ by 
   setting $\alpha'_{e}=\ad(\delta_{e}^{-1})\after\alpha_{e}$, for each $e\in \orE$.
\ede

\ble\label{lem:all pointings are created equally}
Any pointing of $(\cC_0,\orG)$ is isomorphic to itself as a graph of groups.
\ele
\pf
In the definition of an inner morphism let all $\phi_i$ and $\phi_{e}$ be the identity maps and let $\delta_{e}$ be as defined in~\ref{dfn:a pointing of a graph of groups}. This defines an inner isomorphism of graphs of groups.
\qed

\bde\label{dfn:isomorphism of pointings}
An isomorphism between pointings $((\cC,\orG),\delta^{(k)})$ of a graph of groups $(\cC,\orG)$ is an inner isomorphism $\phi$ of $(\cC,\orG)$ such that there exist $a_i\in A_i$ and $a_{e}\in A_{e}$ so that $a_{e}= a_{\bar e}$ and with $\phi_i=\ad(a_i)$ and $\phi_{e}=\ad(a_{e})$ for each $i\in I, e\in \orE$ and we have 
\beq\label{eqn:iso pointing}
\delta_{e}^{(1)}\alpha_{e}(a_{e})= a_{d_0(e)}\delta_{e}^{(2)}.
\eeq
We will say that $\{a_{e}, a_i\mid i\in I, e\in \orE\}$ induces the isomorphism.
\ede

Note that {\em any} choice of $a_e\in A_e$ and $a_i\in A_i$ gives rise to an inner isomorphism of $(\cC,\orG)$. Hence 
Condition~(\ref{eqn:iso pointing}) is the only condition to be checked when verifying if such a collection induces an isomorphism of pointings.

\bde\label{dfn:D}
Consider the reference amalgam $\cA_0=\{G_i,G_{e},\psi_{e}\}$.  Recall that every amalgam of type $\cA_0$ has the same target subgroups $\bar{G}_i$. 
Define a graph of groups as follows. For each $e\in \orE$ we take $A_{d_0(e)}=\Aut(G_{d_0(e)})$, $A_{e}=\Aut_{G_{e}}(\bar G_{d_0(e)}, \bar{G}_{d_1(e)})$ and $\alpha_{e}=\ad(\psi_{e})$. The resulting graph of groups will be denoted by $(\cC_0,\orG)$ and will be called the reference graph.
\ede

\bde\label{dfn:pointings and amalgams}
If $\cA=\{G_i,G_{e}, \varphi_{e}\}$ is an amalgam of type $\cA_0$, we define an associated pointing $((\cC_0,\orG),\delta^\cA)$ as follows.
 If $e\in\orE$ then $\delta^\cA_{e}=\varphi_{e}^{-1}\after\psi_{e}$ where of course by $\varphi_{e}^{-1}$ we mean the inverse of the map $\varphi_{e}\colon G_{d_0(e)}\to \bar G_{d_0(e)}$.
 
 Conversely if $(\cC_0, \orG),\delta)$ is a pointing of $\cC_0$ then we define an amalgam $\cA=\{G_i,G_e,\varphi_{e}\}$ of type $\cA_0$ via  $\varphi_{e}=\psi_{e}\delta_e^{-1}$. 
\ede
 Clearly the two constructions in~Definition~\ref{dfn:pointings and amalgams} are inverse to one another. The questions is of course whether they preserve isomorphism classes.

 \begin{mainthm}\label{thm:pointings and amalgams}
The correspondence $\cA\leftrightarrow ((\cC_0,\orG),\delta^\cA)$ above yields a bijection between isomorphism classes of amalgams and isomorphism classes of pointings of $\cC_0$.
 \end{mainthm}
\pf 
Suppose  $\{a_{e}, a_i\mid i\in I, e\in \orE\}$ induces an isomorphism of pointings $((\cC_0, \orG),\delta')\to((\cC_0, \orG),\delta)$.  We now define a map $\chi$ from $\cA=\{G_i,G_{e}, \varphi_{e}\}$  to $\cA'=\{G_i,G_e, \varphi'_{e}\}$
 by setting $\chi_{e}=a_{e}\in A_e$ for all $e\in E$ and $\chi_i=a_i\in A_i$ for all $i\in I$.  
We claim that $\chi$ is an isomorphism. Indeed, recalling that $ \delta'_{e}\alpha_{e}(\chi_{e})= \chi_{d_0(e)}\delta_{e}$ and $\alpha_e=\ad(\psi_{e})$, we find
$$ \chi_{e}\varphi_{e}= \chi_{e}\psi_{e}\delta_{e}^{-1}=\psi_{e}\alpha_{e}(\chi_e)\delta_{e}^{-1}=\psi_{e}(\delta_{e}')^{-1}\chi_{d_0(e)}=\varphi_{e}'\chi_{d_0(e)}.$$
Conversely if $\chi$ is an isomorphism between $\cA=\{G_i,G_{i,j}, \varphi_{e}\}$ and $\cA'=\{G_i,G_{i,j}, \varphi'_{e}\}$, then the corresponding $a_{e}=\chi_{e}$ and $a_i=\chi_i$ induce isomorphisms between the pointings of the graph. Indeed:
$$\delta'_{e}\alpha_{e}(\chi_{e})=((\varphi'_{e})^{-1}\psi_{e})(\psi_{e}^{-1}\phi_{e}\psi_{e})=(\varphi'_{e})^{-1}\chi_{e}\psi_{e} = \chi_i\varphi_{e}^{-1}\psi_{e}=\chi_i \delta_{e}. $$
\qed

\subsection{Applications}
\paragraph{Goldschmidt's Lemma}

We can now see that that the classical Goldschmidt's Lemma  (see 2.7 of \cite{Gol80} ) is a special case of the above. Recall that an amalgam in the sense of \cite{Gol80} is a pair of group monomorphisms $P_{1} \stackrel{\varphi_{1}}{\leftarrow} B \stackrel{\varphi_{2}}{\rightarrow}   P_{2}$. Consider an amalgam $P_{1} \stackrel{\psi_{1}}{\leftarrow} B \stackrel{\psi_{2}}{\rightarrow}   P_{2}$ and let $A_{i}$ be the subgroup of $\Aut(P_{i})$ that leaves $\im \psi_{i}$ invariant and define $\alpha_{i}\colon A_{i}\to \Aut(B)$ be $\ad(\psi_{i})$ and let $\bar A_{1}= \im \alpha_{i}$. The terminology of loc.~cit.\ agrees with ours in the sense that an amalgam of type $\cA_0=\{P_1,B,P_2,\psi_1,\psi_2\}$ is simply an amalgam where all groups $P_1$, $B$, $P_2$, $\im(\psi_1)$, and $\im(\psi_2)$ are the same, but the connecting monomorphisms $\varphi_1$ and $\varphi_2$ might be different from $\psi_1$ and $\psi_2$.

\bco[Goldschmidt's Lemma] \label{lem:Goldschmidt} There is a 1-1 correspondence between isomorphism classes of amalgams of type $\cA_0$ and double cosets of $\bar A_{1}$, and $\bar A_{2}$ in $\Aut(B)$.
\eco
\pf We can view the amalgamated product as an amalgam associated to a double loop $\orG = (\{a\}, \{e_1, e_2\})$ with $G_{e_i}=P_i$ and $G_a=B$. The associated graph of groups $(\cC_0,\orG)$ is given on the groups $A_{e_{i}}=A_{i}$, $A_{0}=\Aut(B)$ by maps $\alpha_i$, as described above.
It follows from Theorem~\ref{thm:pointings and amalgams} that the isomorphism classes of amalgams of type $\cA_0$ correspond to the isomorphism classes of pointings of the graph of groups $(\cC_0,\orG)$. Given a pointing $((\cC_{0},\orG),\{\delta_{e_{1}}, \delta_{e_{2}}\})$, one immediately notes that the set $\{a_{e_{i}}=\id, a_{0}=\delta_{e_{1}}\mid i=1,2\}$ induces an isomorphism between $((\cC_{0},\orG),\{\delta_{e_{1}}, \delta_{e_{2}}\})$ and $((\cC_{0},\orG),\{\id_{e_{1}},\delta_{e_{1}}^{-1} \delta_{e_{2}}\})$ in the sense of Definition~\ref{dfn:isomorphism of pointings}.
Thus, it suffices to  consider the isomorphism classes of pointings of type $((\cC_{0},\orG),\{\id_{e_{1}}, \delta_{e_{2}}\})$. Moreover two pointings  $((\cC_{0},\orG),\{\id_{e_{1}}, \delta^{(k)}_{e_{2}}\})$ ($k=1,2$) are isomorphic if and only if there exist $a_{i}\in A_{i}$ so that $\delta^{(1)}_{e_{2}}\alpha_{e_{2}}(a_{e_{2}})=\alpha_{e_{1}}(a_{e_{1}})\delta_{e_{2}}^{(2)}$.  Since $\im \alpha_i=\bar A_i$, this happens if and only if $\delta^{(k)}_{e_{2}}$ are in the same $\bar A_{1}, \bar A_{2}$ double coset.\qed

\paragraph{Triangle graphs}
We now treat the special case where $\orG$ is a triangle.
This seems arbitrary but it is a very common example. For instance, the amalgam of maximal parabolics  in a group acting on a rank three geometry has a triangle as its diagram. 
To fix notation, let $I=\{1,2,3\}$ and $E=\{\{1,2\}, \{2,3\}, \{3,1\}\}$,  so we have three edge groups $A_{i,i+1}$ and 3 vertex groups $A_{i}$ (in the rest of the section we will abuse the notations and consider the indices modulo 3) and the respective $\alpha_{i,i+1}\colon A_{i,i+1}\to A_i$ and $\alpha_{i+1,i}\colon A_{i,i+1}\to A_{i+1}$.
Consider the group $A=A_{1}\times A_{2}\times A_{3}$, the inclusion maps $i_{k}\colon A_{k}\to A$ and projection maps  $\pi_{k }\colon A \to A_{k}$. Moreover for each $k$ there is a map $\Phi_{k}\colon A_{k,k+1}\to A$ given by $x\mapsto i_{k}(\alpha_{k,k+1}(x))\cdot i_{k+1}(\alpha_{k+1,k}(x)) $, (for example $\Phi_{1}(x)=(\alpha_{1,2}(x),\alpha_{2,1}(x), 1)$). Consider $H_{k}=\im\Phi_{k}$. Define an equivalence relation $\sim$ on $A$ as follows 
$$(a_{1},a_{2},a_{3})\sim (a_{1}',a_{2}',a_{3}') \iff  \forall k=1,2,3, \exists x_{k}\in H_{k}, \mbox{ such that } a'_{k}=\pi_{k}(x_{k}^{-1}i_{k}(a_{k})x_{k-1})$$
\bco
The equivalence classes for $\sim$ are in bijection to isomorphism classes of pointings of $\cC_{0}$.
\eco
\pf
First note that if $\{\delta_{i,i+1}, \delta_{i+1,i}\mid i=1,2,3\}$ give a pointing of $\cC_{0}$ then the collection $\{a_{i,i+1}=1, a_{i}=\delta_{i,i+1}\mid i = 1,2,3\}$ induces an isomorphism between the pointing $((\cC_{0},\orG),\delta)$ and the pointing $((\cC_{0},\orG),\delta')$ where $\delta'_{i,i+1}=1$ and $\delta'_{i+1,i}=\delta_{i+1,i-1}^{-1}\delta_{i+1,i}$. This means that any pointing is isomorphic to a pointing in which $\delta_{i,i+1} =1$ for all $i=1,2,3$. We call such a pointing {\em positively normalized}.

We now describe a bijection $\Xi$ between isomorphism classes of pointings $((\cC_{0},\orG),\delta)$ and classes in $A/\!\!\sim$.
To the isomorphism class of the pointing given by $\{\delta_{i,i+1},\delta_{i+1,i}\mid i=1,2,3\}$, we associate the equivalence class 
$[(\delta_{1,2}^{-1}\delta_{1,3},\delta_{2,3}^{-1}\delta_{2,1},\delta_{3,1}^{-1} \delta_{3,2})]$. Conversely, to each equivalence class $[(a_{1},a_{2},a_{3})]$ we associate the isomorphism class containing the normalized pointing given by $\{\delta_{i,i+1}=1,\delta_{i+1,i}=a_{i+1}\mid i=1,2,3\}$. 

We now show that these maps are well-defined and each other's inverse.
First we make the following observation. Suppose we have two positively normalized pointings $\delta$ and $\delta'$ and a collection $\{a_{i,i+1},a_{i+1,i},a_i\mid i=1,2,3\}$ inducing an isomorphism $((\cC,\orG),\delta)\to ((\cC,\orG),\delta')$.
Then, since for directed edges $e=(i,i+1)$ those pointings have $\delta_{i,i+1}=\delta_{i,i+1}'=1$, Condition~(\ref{eqn:iso pointing}) turns into 
$$\begin{array}{rl}
\alpha_{i,i+1}(a_{i,i+1})&=a_{i},\\
\delta_{i+1,i}\alpha_{i+1,i}(a_{i+1,i})&=a_{i+1}\delta'_{i+1,i},\\
\end{array}\mbox{ for }i=1,2,3.$$  
Thus, the existence of the set $\{a_{i,i+1},a_{i+1,i},a_i\mid i=1,2,3\}$ is equivalent to the existence of a set $\{a_{i,i+1},a_{i+1,i}\mid i=1,2,3\}$ satisfying
$$\delta'_{i+1,i}=\alpha_{i+1,i-1}(a_{i+1,i-1})^{-1}\delta_{i+1,i}\alpha_{i+1,i}(a_{i+1,i})\mbox{ for }i=1,2,3,$$  
On the other hand, consider $(a_{1},a_{2},a_{3})\sim (a_{1}',a_{2}',a_{3}')$. This is equivalent to the fact that there exist $x_{k}=\Phi_{k}(a_{k,k+1})$ ($k=1,2,3$) so that
 $$
 \begin{array}{@{}l@{}}
 a'_{1}=\pi_{1}(x_{1}^{-1}i_{1}(a_{1})x_{3})=\pi_{1}(\alpha_{1,2}(a_{1,2}^{-1})a_{1}\alpha_{1,3}(a_{3,1}),\alpha_{2,1}(a_{1,2}^{-1}), \alpha_{3,1}(a_{3,1}))=\alpha_{1,2}(a_{1,2}^{-1})a_{1}\alpha_{1,3}(a_{3,1})
 \\
 a'_{2}=\pi_{1}(x_{2}^{-1}i_{2}(a_{2})x_{1})=\pi_{2}(\alpha_{1,2}(a_{1,2}), \alpha_{2,3}(a_{2,3}^{-1})a_{2}\alpha_{2,1}(a_{1,2}),\alpha_{3,2}(a_{2,3}^{-1}))=\alpha_{2,3}(a_{2,3}^{-1})a_{2}\alpha_{2,1}(a_{1,2})
 \\
 a'_{3}=\pi_{3}(x_{3}^{-1}i_{3}(a_{3})x_{2})=\pi_{3}(\alpha_{1,3}(a_{3,1}^{-1}),\alpha_{2,3}(a_{2,3}), \alpha_{3,1}(a_{3,1}^{-1})a_{3}\alpha_{3,2}(a_{2,3}))= \alpha_{3,1}(a_{3,1}^{-1})a_{3}\alpha_{3,2}(a_{2,3})
 \end{array}$$
Thus, we see that the map $\Xi$ is a bijection between isomorphism classes of pointings and equivalence classes in $A$.
\qed

\section{The fundamental group}\label{sec:fundamental group}
In this section we shall restrict to a special kind of graph of groups.

\bde\label{dfn:rigid graph of groups}
In the notation developed in Section~\ref{sec:Amalgams and graphs of groups}, we say that the graph of groups $(\cC,\orG)$ is {\em rigid} if the following condition is satisfied. 
\begin{enumerate}
\item[(Ri)]For any edge $e$, there exists a subgroup $\tilde A_{e}\le A_{e}=\Aut_{G_{e}}(\bar G_{d_0(e)}, \bar G_{d_{1}(e)})$ so that $\alpha_{e}|_{\tilde A_{e}}$  and $\alpha_{\bar e}|_{\tilde A_{\bar e}}$ are bijections onto the images of $\alpha_e$ and $\alpha_{\bar e}$ respectively.
\end{enumerate}
In other words, we require $\ker(\alpha_e)$ to have a complement $\tilde A_e$ for each $e\in \orE$.
We then set $\tilde \cC_{0}=\{A_{i }, \tilde A_{e}, \alpha_{e}\}$.
We call the amalgam $\cA$ rigid if the graph of groups associated to it as in Theorem~\ref{thm:pointings and amalgams}
 is rigid.
\ede

\begin{mainthm}
If the amalgam $\cA$ is rigid then there is a 1-1 correspondence between between isomorphism classes of amalgams of type $\cA$ and isomorphism classes of pointings of $\tilde\cC_0$.

\end{mainthm}

\pf
In view of Theorem~\ref{thm:pointings and amalgams}, we must show that there is a 1-1 correspondence between isomorphism classes of pointings of $\cC_0$ and those of $\tilde \cC_{0}$. The two graphs of groups are described by the same amalgam type and we will consider the map $((\cC_{0},\orG),\delta)\mapsto ((\tilde \cC_{0},\orG),\delta)$ that keeps the same elements of the pointings. Of course if two pointings of $(\tilde \cC_{0},\orG)$ are isomorphic then the corresponding pointings of $\cC_{0}$ are, using the same elements of $\tilde A_{e}\le A_{e}$. 

The converse is also true. Indeed if $\{a_{e}, a_i\mid i\in I, e\in \orE\}$ defines an isomorphism between two pointings of $\cC_{0}$, then  for each $e$ we can pick $\tilde a_{e}\in \tilde A_{e}$ so that  $\alpha_{e}(a_{e})=\alpha_{e}(\tilde a_{e})$. The system $\{\tilde a_{e}, a_i\mid i\in I, e\in \orE\}$induces an isomorphism between the corresponding pointings of $\tilde C_{0}$.
\qed
\mn
Note that $\tilde \cC_{0}$ is a graph of groups in the more restricted sense of \cite{Bas93} and so we can use Definitions~\ref{def:path group}~and~\ref{def:fundamental group} introduced there.

\bde\label{def:path group}
For a given graph of groups $(\cC, \orG)$ we define the {\em path group} as follows
$$ \pi(\cC)=((*_{i\in I} A_i )* F(\orE))/ R$$
where $F(\orE)$ is  the free group on the set $\orE$, $*$ denotes free product and $R$ is the following set of relations
$$ \mbox{for any }  e=(i,j) \in \orE, a \in  A_e , e\bar e=\id \mbox{ and } e\cdot \alpha_{\bar{e}}(a)\cdot\bar e =\alpha_e(a) $$
\ede

\bde\label{def:fundamental group}
Given a graph of groups $(\cC, \orG)$, a {\em path of length $n$} in $\cC$ is a sequence $a_1 e_1 a_2 \cdots e_n a_n$ where $e_1, \ldots, e_n$ are  edges  in $\orG$ with $d_{0}(e_{k})=i_{k}$ and $a_k \in A_{i_k}$.
The sequence $e_1, \ldots, e_n$ is called an {\em edge path}.  A path is {\em reduced} if it has no returns (i.e.\ $e_{i+1}\ne \bar e_i$). This defines an element of   $\pi(\cC)$  by  setting $|\gamma|=a_1\cdot e_1 \cdot a_2 \cdots e_n \cdot a_n $. We denote by $\pi[i,j]=\{|\gamma| \mid \gamma \mbox{ a path on } \cC, \gamma(1)=i, \gamma(2)=j \}$ respectively $\pi(\cC, i)=\pi[i,i]$. Concatenation induces a group operation on $\pi(\cC, i)$ and we call this group {\em the fundamental group of $\cC$ with base point $i$}.    
\ede

\bde\label{lem: fund group of a pointing}
If $((\cC_0, \orG), \delta)$ is a pointing of the graph of groups $(\cC_0, \orG)$ then any  edge path  $\gamma=e_1 \cdots e_n$ in $\orG$ gives rise to a  path in $\cC$ via $\gamma \mapsto \gamma_{\delta}=\delta_{e_1} e_1 \delta_{\bar e_1}^{-1}\delta_{e_2} \cdots e_{n-1}\delta_{\bar e_{n-1}}^{-1}\delta_{e_n}e_n \delta_{\bar e_n}^{-1}$. Fixing a base point $a$, the map $\gamma\mapsto |\gamma_\delta|$ restricts to a monomorphism $\Phi\colon\pi_1( \Gamma, a) \to \pi(\cC_0,a)$, where $\pi_1(\Gamma,a)$ denotes the fundamental group of $\Gamma$ with base point $a$. The image of this map is called {\em the fundamental group of the pointing} and denoted by $\pi(\cC,a,\delta)$. 
\ede

\bth\label{thm:pointings and fundamental group}
Two pointings of $\tilde \cC_{0}$ are isomorphic if and only if they have the same fundamental group.
\eth
\pf
The proof is identical to that of Theorem 3.21 of \cite{BloHof09a}.
\qed
\mn
\Note One should be careful not to misread Theorem~\ref{thm:pointings and fundamental group} as saying that pointings with {\em isomorphic} fundamental groups are isomorphic. This is not the case since all such fundamental groups are isomorphic to the fundamental group $\pi_1(\Gamma,a)$ and not all pointings are isomorphic.

\subsection{Applications}
Consider the reference amalgam $\cA_0=\{G_i,G_{e},\psi_{e}\}$.  Recall that every amalgam of type $\cA_0$ has the same target subgroups $\bar{G}_i$.
For each $e\in \orE$ define 
  $\bar{D}_e=N_{G_e}(\bar{G}_{d_1(e)})\cap \bar{G}_{d_0(e)}$
   and $D_e=\psi_{e}^{-1}(\bar{D}_e)$.
Let us consider a non-collapsing amalgam $\cA=\{G_i,G_{e},\varphi_{e}\}$ of type $\cA_0$ with graph $\orG$. 
In~\cite{BeSh2004,BloHof09a,Dun2005} it is shown that the fact that $\orG$  has no triangles forces the amalgam to have the following property:
\begin{itemize}
\item[(D1)] for any two edges $e,f$ such that 
 $d_0(e)=d_0(f)$ we have $\varphi_{e}^{-1}(\bar{D}_e)=\varphi_{f}^{-1}(\bar{D}_f)$;
 \end{itemize}
In~\cite{Kee2009} the failure of property (D1) is used to show that any amalgam of mixed Phan-Curtis-Tits type over  a graph without triangles collapses.
 
\mn
As a consequence of (D1) one can show that $\cA$ is isomorphic to an amalgam with the following property 
\begin{itemize}
\item[(D2)] $D_e=\varphi_{e}^{-1}(\bar{D}_e)$ for all $e\in \orE$.
\end{itemize}
Thus, in order to classify amalgams of type $\cA_0$ we can restrict to those with properties (D1) and (D2). We can then also write   $D_{d_0(e)}=D_e$, for all $e\in \orE$, without ambiguity.

For any $e\in \orE$ we now have 
 $\delta_e^\cA=\varphi_e^{-1}\after \psi_e\in \Aut_{G_{d_0(e)}}(D_{d_0(e)})$. Hence by Theorem~\ref{thm:pointings and amalgams}, isomorphism classes of such amalgams are in bijection with pointings of the graph of groups 
  $(\cC_0,\orG)=\{A_i,A_e,\alpha_e\}$, where $A_e$ is as before, but $A_i=\Aut_{G_i}(D_i)$.
In the cases of loc.~cit.~this graph of groups is in fact rigid.

\end{document}